\documentclass[11pt]{article}
\usepackage{amssymb,amsfonts,amsmath,latexsym,epsf,tikz,url}
\newtheorem{theorem}{Theorem}[section]

\newtheorem{corollary}[theorem]{Corollary}

\newcommand{\proof}{\noindent{\bf Proof.\ }}
\newcommand{\qed}{\hfill $\square$\medskip}

\begin{document}

\title{ Total dominator chromatic number of some operations on a graph}

\author{Nima Ghanbari and Saeid Alikhani$^{}$\footnote{Corresponding author} }

\date {\today}

\maketitle

\begin{center}

   Department of Mathematics, Yazd University, 89195-741, Yazd, Iran\\
{\tt n.ghanbari.math@gmail.com, alikhani@yazd.ac.ir}\\

\end{center}


\begin{abstract}
Let $G$ be a simple graph. A total dominator coloring of $G$ is a proper coloring of the vertices of $G$ in which each vertex of the graph is adjacent to every vertex of some color class.
The total dominator chromatic number $\chi_d^t(G)$ of $G$ is the minimum number of colors
among all total dominator coloring of $G$. In this paper, we  examine the effects on $\chi_d^t(G)$ when $G$ is modified by operations on vertex and  edge of $G$.  
\end{abstract}


\medskip
\noindent{\bf AMS Subj.\ Class.:} 05C15, 05C69

\section{Introduction}

In this paper, we consider  simple finite graphs, without directed, multiple, or weighted edges, and without self-loops. Let $G=(V,E)$ be such a graph and $k \in \mathbb{N}$. A mapping $f : V (G)\longrightarrow \{1, 2,...,k\}$ is
called a $k$-proper  coloring of $G$ if $f(u) \neq f(v)$ whenever the vertices $u$ and $v$ are adjacent
in $G$. A color class of this coloring is a set consisting of all those vertices
assigned the same color. If $f$ is a proper coloring of $G$ with the coloring classes $V_1, V_2,..., V_{k}$ such
that every vertex in $V_i$ has color $i$, then sometimes write simply $f = (V_1,V_2,...,V_{k})$.  The chromatic number $\chi(G)$ of $G$ is
the minimum number of colors needed in a proper coloring of a graph. The chromatic number is perhaps the most studied of all graph theoretic parameters.
A dominator coloring of $G$ is a proper coloring of $G$ such that every vertex
of $G$ dominates all vertices of at least one color class (possibly its own class), i.e., every vertex of $G$ is adjacent to all vertices of at least one color class.  The dominator
chromatic number $\chi_d(G)$ of $G$ is the minimum number  of color classes in a dominator coloring of $G$.
Kazemi \cite{Adel,Adel2} studied the total dominator coloring, abbreviated TD-coloring. Let $G$ be a graph with no
isolated vertex, the total dominator coloring  is a proper coloring of $G$ in which each vertex of the graph is adjacent
to every vertex of some (other) color class. The total dominator chromatic number, abbreviated TD-chromatic number, $\chi_d^t(G)$ of $G$ is the minimum number of color classes in a TD-coloring of $G$. The TD-chromatic number of a graph is related to its total domination
number. Recall that a total dominating set of $G$ is a set $S\subseteq V(G)$ such
that every vertex in $V(G)$ is adjacent to at least one vertex in $S$ and the total domination
number of $G$, denoted by $\gamma_t(G)$, is the minimum cardinality of a total dominating set of $G$. A total dominating set of $G$ of cardinality $\gamma_t(G)$ is called a $\gamma_t(G)$-set. The literature on the subject on total domination in graphs
has been surveyed and detailed in the  book \cite{Henningbook}. It has proves that the computation of the TD-chromatic number is NP-complete (\cite{Adel}).
The TD-chromatic number of some graphs, such as paths, cycles, wheels and the complement of paths and cycles has computed in
\cite{Adel}. Also Henning in \cite{GCOM} established the  lower and upper bounds on the TD-chromatic number
of a graph in terms of its total domination number. He has shown that, for  every
graph $G$ with no isolated vertex satisfies $\gamma_t(G) \leq \chi_d^t (G)\leq \gamma_t(G) + \chi(G)$.
The properties of TD-colorings in trees has studied in \cite{GCOM,Adel}. Trees $T$ with $\gamma_t(T) =\chi_d^t(T)$ has characterized
in \cite{GCOM}. In \cite{Alikhani} considered graphs with specific construction and study their TD-chromatic number.  
The join $G = G_1 + G_2$ of two graph $G_1$ and $G_2$ with disjoint vertex sets $V_1$ and $V_2$ and
edge sets $E_1$ and $E_2$ is the graph union $G_1\cup G_2$ together with all the edges joining $V_1$ and
$V_2$. For two graphs $G = (V,E)$ and $H=(W,F)$, the corona $G\circ H$ is the graph arising from the
disjoint union of $G$ with $| V |$ copies of $H$, by adding edges between
the $i$th vertex of $G$ and all vertices of $i$th copy of $H$.
In the study of TD-chromatic number of graphs,  this naturally raises the question:  What happens to the TD-chromatic number, when we consider some operations on the vertices and  the edges of a graph? 
In this paper we would like to answer to this question. 
\medskip

In the next section, examine the effects on $\chi_d^t(G)$ when $G$ is modified by deleting a vertex or deleting an edge. In Section 3, we   study the effects on $\chi_d^t(G)$,  when $G$ is modified by contracting a vertex and contracting an  edge. Also we consider another obtained graph  by operation on a vertex $v$ denoted by $G\odot v$ which is a graph obtained from $G$ by the removal of all edges between any pair of neighbors of $v$ in Section 3 and study $\chi_d^t(G\odot v)$.


\section{Vertex and edge removal}

The graph $G-v$ is a graph that is made  by deleting the vertex $v$ and all edges connected to $v$ from the graph $G$ and the graph $G-e$ is a graph that obtained from $G$ by simply removing the  edge $e$. Our main results in this section are  in obtaining  a  bound for TD-chromatic number  of $G-v$ and $G-e$. To do this, we need to consider some preliminaries.
\begin{theorem}{\rm(\cite{Adel})}\label{CnPn}
	\item[(i)]
	Let $P_n$ be a path of order $n\geq 2$. Then
	\[
	\chi_d^t(P_n)=\left\{
	\begin{array}{lr}
	{\displaystyle
		2\lceil\frac{n}{3}\rceil-1}&
	\quad\mbox{if $n\equiv 1$ $(mod\,3)$,}\\[15pt]
	{\displaystyle
		2\lceil\frac{n}{3}\rceil}&
	\quad\mbox{otherwise.}
	\end{array}
	\right.
	\]
	\item[(ii)] Let $C_n$ be a cycle of order $n\geq 3$. Then
	\[
	\chi_d^t(C_n)=\left\{
	\begin{array}{lr}
	{\displaystyle
		2} &
	\quad\mbox{if $n=4$}\\[15pt]
	{\displaystyle
		4\lfloor\frac{n}{6}\rfloor+r}&
	\quad\mbox{if $n\neq 4$, $n\equiv r$ $(mod\,6)$, $r=0,1,2,4$,}\\[15pt]
	{\displaystyle
		4\lfloor\frac{n}{6}\rfloor+r-1}&
	\quad\mbox{if $n\equiv r$ $(mod\,6)$, $r=3,5$.}
	\end{array}
	\right.
	\]
\end{theorem}

The following theorem gives an upper bound and a lower bound  for  $\chi_d^t(G-e)$.  
\begin{theorem}\label{G-e1} 
Let $G$ be a connected graph, and $e=vw\in E(G)$ is not a bridge of $G$. Then we have:
	$$\chi_d^t(G)-1\leq \chi_d^t(G-e)\leq \chi_d^t(G)+2 .$$
\end{theorem}

\proof 
First we prove the left inequality. We shall present  a TD-coloring for $G-e$. If we add the edge $e$ to $G-e$, then we have two cases. If two vertices  $v$ and $w$ have the same color in the TD-coloring of $G-e$, then in this case we add a new color, like $i$, to one of them. Since every vertex use the old class for TD-coloring then this is  a TD-coloring for $G$. So we have $\chi_d^t(G)\leq \chi_d^t(G-e) +1.$
If two vertices  $v$ and $w$ do not have the same color in the TD-coloring of $G-e$, then  the TD-coloring of $G-e$ can be a TD-coloring for $G$. So $\chi_d^t(G)\leq \chi_d^t(G-e)$ and therefore we have $\chi_d^t(G)-1\leq \chi_d^t(G-e).$

\medskip
Now we prove $\chi_d^t(G-e)\leq \chi_d^t(G)+2$. Suppose that  the vertex $v$ has color $i$ and $w$ has color $j$. We have the following cases:

Case 1) The vertex $v$ does not use the color class $j$ and $w$ does not use the color class $i$ in the TD-coloring of $G$. So the TD-coloring of $G$ gives a TD-coloring of $G-e$ and in this case $\chi_d^t(G-e)=\chi_d^t(G)$.

Case 2) The vertex $v$  uses the color class $j$ but $w$ does not use the color class $i$ in the TD-coloring of $G$. Since $v$ used the color class $j$ for the TD-coloring then we have two cases:
\begin{enumerate}
	\item[(i)]
	If $v$ has some adjacent vertices which have color $j$, then we give the new color $l$ to all of these vertices and this coloring is a TD-coloring for $G-e$.
	\item[(ii)]	
	If any other vertex does not have color $j$, since $G-e$ is a connected graph, then exists vertex $s$ which is adjacent to $v$. Now we give to $s$ the new color $l$ and this coloring is a TD-coloring for $G-e$.
\end{enumerate}
So for this case, we have $\chi_d^t(G-e)= \chi_d^t(G)+1$.

Case 3)  The vertex $v$  uses the color class $j$ and $w$  uses the color class $i$ in the TD-coloring of $G$. We have three cases:
\begin{enumerate}
	\item[(i)]
	There are some vertices which are adjacent to $v$ and have color $j$. Then we color all of them with color $l$. And there are some vertices which are adjacent to $w$ and have color $i$. We color all of them with color $k$. So this is a TD-coloring for $G-e$.
	\item[(ii)]
	Any other vertex does not have color $j$. Then we do the same as Case 2 (ii) and  there are some vertices which are adjacent to $w$ and have color $i$. Then we do the same as Case 3 (i).
	\item[(iii)]
	Any other vertex does not have colors $i$ and $j$. Then we do the same as Case 2 (ii) and use two new colors $l$ and $k$.
\end{enumerate}
So we have $\chi_d^t(G-e)\leq \chi_d^t(G)+2$. 
\qed

\begin{figure}
	\begin{center}
		\includegraphics[width=3.5in]{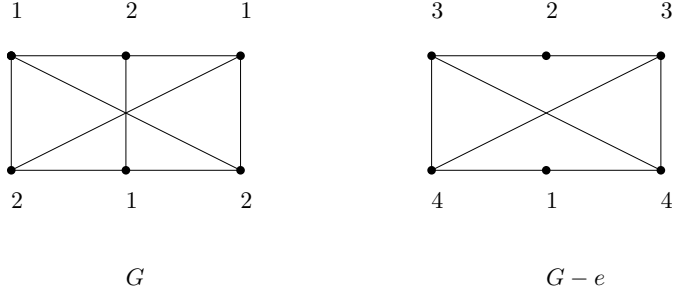}
		\caption{TD-coloring of $G$ and $G-e$.}
		\label{GGe}
	\end{center}
\end{figure}

\medskip 
\noindent{\bf Remark 1.} The lower bound of $\chi_d^t(G-e)$  in Theorem \ref{G-e1} is sharp. It suffices to consider complete graph $K_3$ as $G$. Also the upper bound is sharp, because as we see in Figure \ref{GGe},  $\chi_d^t(G)=2$ and $\chi_d^t(G-e)=4$.

\medskip 
Now we consider the graph $G-v$, and present a lower bound and an upper bound for the TD-chromatic number of $G-v$. 

\begin{theorem}\label{G-v1}
Let $G$ be a connected graph, and  $v\in V(G)$ is not a cut vertex of $G$.  Then we have:
	$$\chi_d^t(G)-2\leq \chi_d^t(G-v) \leq \chi_d^t(G)+deg(v)-1 .$$
\end{theorem}

\proof
First we prove $\chi_d^t(G)-2\leq \chi_d^t(G-v)$. We shall present  a TD-coloring for $G-v$. If we add vertex $v$ and all the corresponding edges to $G-v$, then it suffices to give the new color $i$ to vertex $v$ and the new color $j$ only to one of the adjacent vertices of  $v$ like $w$ and do not change all the other colors. Since every vertices except $v$ and $w$ use the old classes for TD-coloring and $v$ uses the color class $j$ and $w$  uses the color class $i$ so we have a TD-coloring of $G$. Therefore we have $\chi_d^t(G)\leq \chi_d^t(G-v) +2 $ and we have the result. 

Now we prove $\chi_d^t(G-v)\leq \chi_d^t(G)+deg(v)-1$. First we give a TD-coloring to $G$. Suppose that the vertex $v$ has the color $i$. So we have the following cases:

Case 1)  There is another vertex with color $i$. In this case every vertex uses the old class for TD-coloring and then this is  a TD-coloring for $G-v$. So $\chi_d^t(G-v)\leq \chi_d^t(G).$

Case 2) There is no other vertex with color $i$. In this case we give the new colors $i,a_1,a_2,\ldots,a_{deg(v)-1}$ to all the adjacent vertices of $v$. Obviously, this is a TD-coloring for $G-v$. Therefore $\chi_d^t(G-v)\leq \chi_d^t(G)+deg(v)-1 .$
\qed

\noindent{\bf Remark 2.} The lower bound in Theorem \ref{G-v1} is sharp. Consider the cycle  $C_{10}$, as $G$. For every $v\in V(C_{10})$ we have $C_{10}-v=P_9$ which is a path graph of order $9$. Then by the Theorem \ref{CnPn} we have $\chi_d^t(C_{10})=8$ and $\chi_d^t(P_9)=6$.

\medskip

\begin{figure}
		\begin{center}
			\includegraphics[width=0.8in]{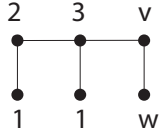}
			\caption{$P_3\circ K_1$}
			\label{P3K}
		\end{center}
	\end{figure}

To obtain more results, we consider the corona of $P_n$ and $C_n$ with $K_1$. The following theorem gives the TD-chromatic number of
these kind of graphs:

\begin{theorem}\label{centi}
	\begin{enumerate}
		\item[(i)] For every $n\geq 2$, $\chi_d^t(P_n\circ K_1)=n+1$.
		\item[(ii)] For every $n\geq 3$, $\chi_d^t(C_n\circ K_1)=n+1$.
	\end{enumerate}
\end{theorem}
\proof
\begin{enumerate}
	\item[(i)]
	We color the $P_n\circ K_1$ with numbers $1,2,...,n+1$, as shown in the Figure \ref{PKCK}. Observe that, we need $n+1$ color for TD-coloring. We shall show that we are not able to have TD-coloring with less colors.
	
	\begin{figure}[h]
		\hspace{1.3cm}
		\begin{minipage}{6.1cm}
			\includegraphics[width=\textwidth]{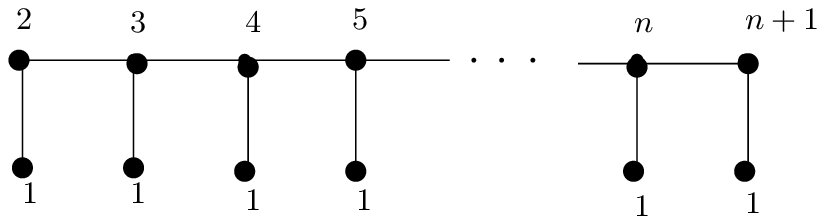}
		\end{minipage}
		\hspace{.7cm}
		\begin{minipage}{6.1cm}
			\includegraphics[width=\textwidth]{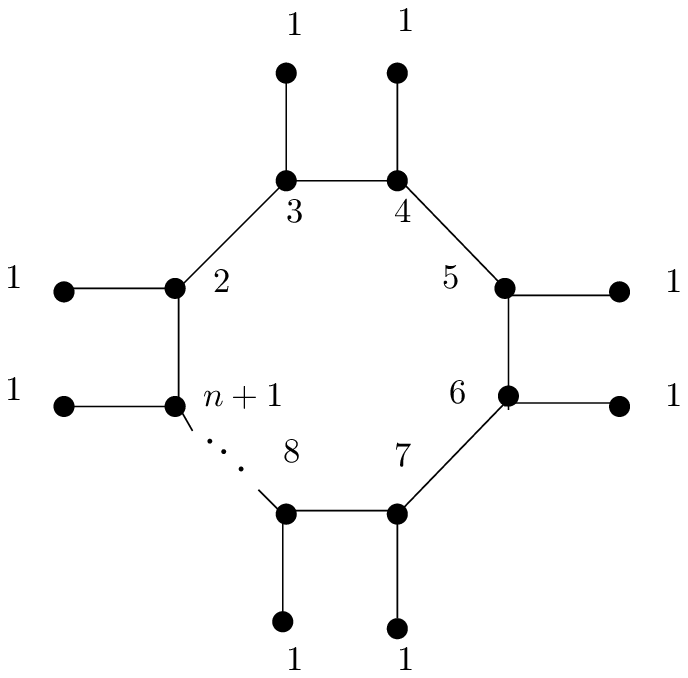}
		\end{minipage}
		\caption{\label{PKCK} Total dominator coloring of $P_n\circ K_1$ and $C_n\circ K_1$, respectively.}
	\end{figure}
	
	Obviously we have $\chi_d^t(P_2\circ K_1)=3$. Now we consider $P_3 \circ K_1$. As we see in Figure \ref{P3K}, we can not give number $1$ to vertex $v$, because there is no number to color vertex $w$.
	Also we can't consider number $2$ for vertex $v$ since the vertex which has color $1$ and is adjacent to vertex with number $2$, is not adjacent with $v$. Since the coloring is proper, we cannot use color $3$ too for this vertex. So we give number $4$ to vertex $v$. Between used colors, we can use only number $1$ for vertex $w$. Therefore $\chi_d^t(P_3\circ K_1)=4$. Similarly, we color $P_i\circ K_1$ from $P_{i-1}\circ K_1$ when $i\geq 3$. Any other kinds of coloring of this graph needs more colors. So we have the result.

\item[(ii)] It is similar to the part (i). \qed
	\end{enumerate} 
	
	\begin{figure}
		\begin{center}
			\includegraphics[width=1.6in]{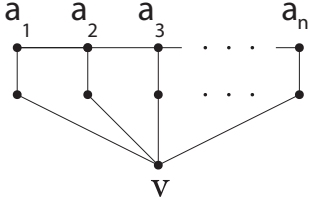}
			\caption{Graph $G$ in the proof of Theorem \ref{large}}
			\label{G1}
		\end{center}
	\end{figure}

We end this section with the following theorem:

\begin{theorem}\label{large}
There is a connected graph $G$, and a vertex  $v\in V(G)$ which is not a cut vertex of $G$ such that 	 $\vert \chi_d^t(G) -  \chi_d^t(G-v) \vert $ can be arbitrarily large.
\end{theorem}

\proof
Consider the graph $G$ in Figure \ref{G1}. We color the vertices $a_1, a_2,\ldots,a_n$ with $\chi_d^t(P_n)$ colors. Then we give the new color $\chi_d^t(P_n)+1$ to all the adjacent vertices of $v$ and $\chi_d^t(P_n)+2$ to $v$. Obviously this is a TD-coloring for $G$. So we have:

\[\chi_d^t(G)=2+
	\chi_d^t(P_n)=\left\{
	\begin{array}{lr}
	{\displaystyle
		2\lceil\frac{n}{3}\rceil+1}&
	\quad\mbox{if $n\equiv 1$ $(mod\,3)$,}\\[15pt]
	{\displaystyle
		2\lceil\frac{n}{3}\rceil+2}&
	\quad\mbox{otherwise.}
	\end{array}
	\right.
	\]

Now by removing the vertex $v$, we have $G-v=P_n\circ K_1$ and by Theorem \ref{centi} we have $\chi_d^t(G-v)=n+1$. So we conclude that $\vert \chi_d^t(G) -  \chi_d^t(G-v) \vert $ can be arbitrarily large.\qed  

\section{Vertex and edge contraction}
 Let $v$ be a vertex in graph $G$. The contraction of $v$ in $G$ denoted by $G/ v$ is the graph obtained by deleting $v$ and putting a clique on the (open) neighbourhood of $v$. Note that this operation does not create parallel edges; if two neighbours of $v$ are already adjacent, then they remain simply adjacent (see \cite{Walsh}).  In a graph $G$, contraction of an edge $e$ with endpoints $u,v$ is the replacement of $u$ and $v$ with a single vertex such that edges incident to the new vertex are the edges other than $e$ that were incident with $u$ or $v$. The resulting graph $G/e$ has one less edge than $G$ (\cite{Bondy}). We denote this graph by $G/e$. 
In this section we  examine the effects on $\chi_d^t(G)$ when $G$ is modified by an edge contraction and vertex contraction. First we consider edge contraction:

\begin{theorem}\label{Goe}
	Let $G$ be a connected graph and $e\in E(G)$. Then we have:
	$$\chi_d^t(G)-2 \leq \chi_d^t(G/e) \leq \chi_d^t(G)+1 .$$
\end{theorem}

\proof
First, we find a TD-coloring for $G$. Suppose that the end points of $e$ are $u$ and $v$. The vertex $u$ has the color $i$ and the vertex $v$ has the color $j$. We give all the used colors in the previous coloring to the vertices $E(G)-\{u,v\}$. Now we give the new color $k$ to $u=v$.  Every vertices on the edges of $E(G)-\{u,v\}$ can uses the previous color class (or even $k$) in this coloring. The vertex $u=v$ uses the color class which used for $u$ or $v$ unless $u$ used the color class $j$ and  $v$ used the color class $i$. In this case, if there is another vertex with color $i$, then $u=v$ uses color class $i$ and if there is another vertex with color $j$, then $u=v$ uses color class $j$. If any other vertex does not have the color $i$ and $j$, then it suffices to give color $i$ to one of the adjacent vertices of $u$ (or $v$) in $G$. Then this is a TD-coloring for $G/ e$. So we have $\chi_d^t(G/ e) \leq \chi_d^t(G)+1$.

To find the lower bound,  we shall give a TD-coloring to $G/e$. We add the removed vertex and all the corresponding edges to $G/ e$ and keep the old coloring for the new graph. Now we consider the endpoints of $e$  and remove the used color. Now add new colors $i$ and $j$ to these vertices. All the vertices of edges in $E(G)-\{u,v\}$ can use the previous color class and $u$ can use color class $j$ and $v$ can use color class $i$. So this is a TD-coloring and we have $\chi_d^t(G)\leq\chi_d^t(G/e)+2. $ Therefore $\chi_d^t(G)-2\leq\chi_d^t(G/e). $
\qed

\noindent{\bf Remark 3.} The bounds in Theorem \ref{Goe} are sharp. For the upper bound consider the cycle  $C_4$ as $G$ and for the lower bound consider cycle $C_5$. 

\begin{corollary}
Suppose that  $G$ is  a connected graph and  $e\in E(G)$ is not a bridge of $G$. We have:
$$\frac{\chi_d^t(G-e)+\chi_d^t(G/e)-3}{2}\leq \chi_d^t(G)\leq\frac{\chi_d^t(G-e)+\chi_d^t(G/e)+3}{2}$$
\end{corollary}

\proof
It follows from Theorems \ref{G-e1} and \ref{Goe}.
\qed  

Now we consider the vertex contraction of graph $G$ and examine the effect on $\chi_d^t(G)$ when $G$ is modified by this operation:  

\begin{theorem}\label{G/v}
	Let $G$ be a connected graph and $v\in V(G)$. Then we have:
	$$\chi_d^t(G)-2 \leq \chi_d^t(G/ v) \leq \chi_d^t(G)+deg(v)-1 .$$
\end{theorem}

\proof
First we present  a TD-coloring for $G$. We remove the vertex $v$ and create $G/v$. We consider  one of the adjacent vertices of $v$ like $u$ and do not  change its color and give the new colors $i,i+1,\ldots,i+deg(v)-1$ to other adjacent vertices  of $v$. Now each vertex which was not adjacent to $v$ can use the previous color class (or if the color class changed, the new color class we give to adjacent vertices of $v$). Therefore we have $\chi_d^t(G/ v) \leq \chi_d^t(G)+deg(v)-1$.

To find the lower bound, at first we shall give a TD-coloring to  $G/v$. We add the vertex $v$, add all the removed edges and remove all the added edges. It suffices to give the vertex  $v$ the new color $i$ and only to one of its adjacent vertices like $w$ the new color class $j$. All the vertices which are not adjacent to $v$ can use the previous color classes. All the adjacent vertices of $v$ can use the color class $i$ and $v$ can use the color class $j$. So we have $\chi_d^t(G) \leq \chi_d^t(G/ v) +2.$ Therefore we have the result.
\qed

\noindent{\bf Remark 4.} The bounds in Theorem \ref{G/v} are sharp. For the upper bound consider the complete bipartite graph $K_{2,4}$ as $G$. We have $\chi_d^t(K_{2,4})=2$. By choosing a vertex which is adjacent to four vertices as $v$, we have $K_{2,4}/v=K_5$ which is the complete graph of order $5$ and $\chi_d^t(K_5)=5$. For the lower bound, we consider  cycle graph $C_5$. For every $v\in V(C_5)$ we have $C_5 /v =C_4$. Now by Theorem \ref{CnPn} we have the result.

\begin{corollary}
Let $G$ be a connected graph. For every $v\in V(G)$ which is not cut vertex of  $G$, we have:

$$\frac{\chi_d^t(G-v)+\chi_d^t(G/ v)}{2}-deg(v)+1\leq \chi_d^t(G)\leq\frac{\chi_d^t(G-v)+\chi_d^t(G/ v)}{2}+2.$$
\end{corollary}

\proof
It follows from Theorems \ref{G-v1} and \ref{G/v}.
\qed  

Here we consider another operation on vertex of a graph $G$ and  examine the effects on $\chi_d^t(G)$ when we do this operation. 
 We denote by $G\odot v$ the graph obtained from $G$ by the removal of all edges between
 any pair of neighbors of $v$, note $v$ is not removed from the graph \cite{alikhani1}. The following theorem gives upper bound and lower bound for $\chi_d^t(G\odot v)$. 

	\begin{figure}
		\begin{center}
			\includegraphics[width=4.15in]{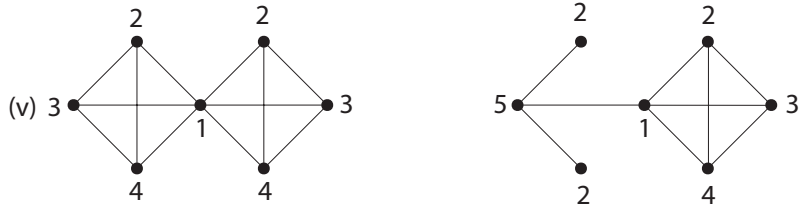}
			\caption{TD-coloring of the graph $G$ and $G\odot v$.}
			\label{G2}
		\end{center}
	\end{figure}

\begin{theorem}\label{godotv}
	Let $G$ be a connected graph and $v\in V(G)$. Then we have:
	$$\chi_d^t(G)-deg(v)+1 \leq \chi_d^t(G\odot v) \leq \chi_d^t(G)+1 .$$
\end{theorem}

\proof
First we prove $\chi_d^t(G\odot v) \leq \chi_d^t(G)+1$. We give a TD-coloring for the  graph $G$. Suppose that the vertex $v$ has the color $i$. We have the following cases:

Case 1) The color  $i$ uses only for the vertex $v$. In this case, adjacent vertices  of the vertex $v$, can use  the color class $i$ and  all the other vertices can use the old color class. So we have $\chi_d^t(G\odot v) \leq \chi_d^t(G)$.

Case 2) The color  $i$ uses  for  another vertex except $v$. In this case, we give the new color $j$ to all of these  vertices (except $v$). This is a TD-coloring for $G\odot v$, because  if a vertex is adjacent to $v$, it can use the color class $i$ and all the other vertices can use old color class and if the old color class changes to $j$ can use $j$ as new color class. So we have $\chi_d^t(G\odot v) \leq \chi_d^t(G)+1$.

Now we prove $\chi_d^t(G)-deg(v)+1 \leq \chi_d^t(G\odot v)$. Consider  the graph $G\odot v$ and shall find a TD-coloring for it. We make $G$ from $G\odot v$ and just change the color of all the adjacent vertices of   $v$ except one of them like $w$ to the new colors $a_1,a_2,\ldots,a_{deg(v)}-1$ and do not change the color of $v$, $w$ and other vertices. This is a TD-coloring for $G$, because $v$ can use the the color class $a_1$. Adjacent vertices of  $v$, can  use the old color class of the TD-coloring of $G\odot v$, and other vertices can use old color class and if the old color classes changes to $a_1$ or $a_2$ or $\ldots$ or $a_{deg(v)-1}$ can use $a_1$ or $a_2$ or $\ldots$ or $a_{deg(v)-1}$ as new color classes. So we have $\chi_d^t(G) \leq \chi_d^t(G\odot v)+deg(v)-1$. Therefore we have the result.
\qed

\noindent{\bf Remark 5.} The bounds in Theorem \ref{godotv} are sharp. For the upper bound consider the graph $G$ in Figure \ref{G2}. It is easy to see that these colorings are TD-coloring.  For the lower bound consider to the complete graph $K_n$ as $G$ ($n\geq 3$). $\chi_d^t(K_n)=n$. Now for every $v\in V(K_n)$, $K_n\odot v$ is the star graph $S_n$ and we have $\chi_d^t(S_n)=2$. By this example we have the following result:

\begin{corollary}
There is a connected graph  $G$  and $v\in V(G)$ such that  $\frac{\chi_d^t(G)}{\chi_d^t(G\odot v)} $ can be arbitrarily large.
\end{corollary}


\begin{thebibliography}{99}
		
\bibitem{alikhani1} S. Alikhani  and E. Deutsch, {\it More on domination polynomial and domination root}, Ars Combin., in press. Available at \texttt{http://arxiv.org/abs/1305.3734v2}.
		
	
   \bibitem{Alikhani} S. Alikhani, N. Ghanbari, {\it  Total dominator chromatic number of specific  graphs}, submitted. Available at \texttt{http://arxiv.org/abs/1511.01652}
   
   \bibitem{Bondy} J.A. Bondy and U.S.R. Murty, {\it Graph Theory with Applications}, American Elsevier, MacMillan, New York, London, (1976).
   
   
	\bibitem{Henningbook} M.A. Henning, A. Yeo, {\it  Total domination in graphs} (Springer Monographs in Mathematics). (2013,
	ISBN: 978-1-4614-6524-9 (Print) 978-1-4614-6525-6 (Online)).
	
	\bibitem{GCOM} M.A. Henning, {\it Total dominator colorings and total domination in graphs}, Graphs Combin. (2015) 31:953-–974.


	\bibitem{Adel} A.P. Kazemi, {\it Total dominator chromatic number of a graph}, Trans. Combin. Vol. 4 No. 2 (2015), pp. 57--68.
	
	
	\bibitem{Adel2} A.P. Kazemi, {\it Total dominator coloring in product graphs}, Utilatas Math., 94 (2014) 329--345. 
	
	 \bibitem{Walsh}	M. Walsh, {\it The hub number of a graph}, Int. J. Math. Comput. Sci., 1 (2006) 117-124.
\end{thebibliography}
\end{document}